\documentclass[12pt,reqno]{amsart}

\usepackage{diagrams,amssymb}
\diagramstyle[height=2em,width=2em, PostScript]
\newarrow{Epi}----{>>}
\newarrow{Mono}>--->
\newarrow{Iso}>---{>>}
\newarrow{Inc}C--->
\newarrow{Mapsto}|--->
\newarrow{Igual}=====
\newarrow{Dashto}{}{dash}{}{dash}>

\newcommand{\Z}{{\mathbb Z}}

\newcommand{\F}{{\mathbb F}}

\newcommand{\pcom}{_{p}^{\wedge}}

\newcommand{\A}{\ifmmode{\mathcal{A}}\else${\mathcal{A}}$\fi}
\newcommand{\K}{\ifmmode{\mathcal{K}}\else${\mathcal{K}}$\fi}
\newcommand{\U}{\ifmmode{\mathcal{U}}\else${\mathcal{U}}$\fi}
\newcommand{\T}{\ifmmode{\mathcal{T}}\else${\mathcal{T}}$\fi}
\newcommand{\FF}{\ifmmode{\mathcal{F}}\else${\mathcal{F}}$\fi}
\newcommand{\LL}{\ifmmode{\mathcal{L}}\else${\mathcal{L}}$\fi}

\newtheorem{Thm}{Theorem}[section]
\newtheorem{Prop}[Thm]{Proposition}
\newtheorem{Cor}[Thm]{Corollary}

\theoremstyle{definition}

\newtheorem{Rmk}[Thm]{Remark}

\theoremstyle{remark}

\title[A spin-off of Serre's theorem] {Relating Postnikov pieces with the Krull
filtration:\\ A spin-off of Serre's theorem}

\author{Nat\`{a}lia Castellana}
\author{Juan A. Crespo}
\author{J\'er\^{o}me Scherer}

\thanks{All three authors are partially supported by MEC grant MTM2004-06686.
The third author is supported by the program Ram\'on y Cajal, MEC,
Spain.}

\begin{document}


\begin{abstract}
We characterize $H$-spaces which are $p$-torsion Postnikov pieces
of finite type by a cohomological property together with a
necessary acyclicity condition. When the mod $p$ cohomology of an
$H$-space is finitely generated as an algebra over the Steenrod
algebra we prove that its homotopy groups behave like those of a
finite complex.
\end{abstract}


\maketitle


\section*{Introduction}
\label{sec intro}
When does cohomological information allow to determine whether or
not a given space is a Postnikov piece? In the $50$'s Serre showed
that a non-trivial $1$-connected finite complex cannot be a
Postnikov piece. He also proved that the same happens for
CW-complexes with finite mod $2$ cohomology \cite{MR0060234}, and
predicted the same behavior at odd primes. After Miller's solution
to Sullivan's conjecture \cite{Miller}, this was proved for spaces
with finite mod $p$ cohomology by McGibbon and Neisendorfer
in~\cite{MR749108}.

The discovery of Lannes' $T$-functor enabled to extend this result
to a larger family of spaces. Indeed, Lannes and Schwartz proved
in~\cite{MR827370} that non-trivial $1$-connected spaces with
locally finite mod~$p$ cohomology cannot be Postnikov pieces.
Finally, in~\cite{MR92b:55004}, Dwyer and Wilkerson showed that,
in fact, this is true for $2$-connected spaces for which the
module of indecomposable elements in the mod $p$ cohomology is
locally finite, including in particular the case where the
cohomology is finitely generated as an algebra.

Observe that the locally finite unstable modules form the $0$th
stage of the Krull-Schwartz filtration $\{\mathcal{U}_n\}$ of the
category $\mathcal{U}$ of unstable modules over the Steenrod
algebra~$\A_p$. This filtration has been introduced in relation
with Kuhn's realizability conjectures, see \cite{MR96i:55027} and
\cite{MR2002k:55043}. Recall that an unstable module $M$ over the
Steenrod algebra lies in $\U_n$ if and only if $\overline T^{n+1}
M = 0$, see~\cite[Theorem~6.2.4]{MR95d:55017}, where $\overline T$
denotes the reduced version of Lannes' $T$ functor.

Thus, the Dwyer-Wilkerson result deals with $2$-connected spaces
$X$ such that $QH^*(X;\F_p)\in \mathcal{U}_0$. In this context we
obtain the following extension for $H$-spaces.

\vskip 4mm

\noindent\textbf{ Theorem~\ref{prop Hdichotomy}.}
{\it Let $X$ be an $(n+2)$-connected $H$-space for some integer $n
\geq 0$ such that $T_V H^* (X; \F_p)$ is of finite type for any
elementary abelian $p$-group $V$. Assume that $QH^*(X; \F_p)$ lies
in $\U_n$. Then either $X$ is contractible, or it has infinitely
many non-trivial homotopy groups with $p$-torsion. In the second
case the iterated loop space $\Omega^{n+1} X$ has infinitely many
non-trivial $k$-invariants.}

\vskip 4mm

Serre's result and its generalizations state conditions on the mod
$p$ cohomology to ensure that a space \emph{is not} a Postnikov
piece. Since mod $p$ cohomology does not detect $q$-primary
information for primes $q \neq p$, our next objective is to give
conditions to ensure that a space is a ($p$-torsion) Postnikov
piece.

It is well-known that $p$-torsion Eilenberg-Mac Lane spaces are
$B\Z/p$-acyclic, that is, their $B\Z/p$-nullification is
contractible (we refer the reader to the book \cite{Dror} for
details about nullification). This implies that $p$-torsion
Postnikov pieces are $B\Z/p$-acyclic as well, so that a first test
to find out if a $p$-torsion space is a Postnikov piece would be
to apply the nullification functor $P_{B\Z/p}$. However , this is
not a sufficient condition as illustrated by the obvious example
of $\prod_{n\geq 1} K(\Z/p, n)$. When dealing with $H$-spaces, we
offer a necessary and sufficient condition in terms of cohomology
and nullification.

\vskip 4mm

\noindent\textbf{ Theorem~\ref{thm Postnikovpieces}.}
{\it Let $X$ be an $H$-space. Then $X$ is a $p$-torsion Postnikov
piece of finite type if and only if $P_{B\Z/p} X$ is contractible
and $H^* (X; \F_p)$ is a finitely generated algebra over~$\A_p$.}

\vskip 4mm

Other examples of $H$-spaces with finitely generated cohomology as
an algebra over the Steenrod algebra are the highly connected
covers of finite $H$-spaces. It follows from Neisendorfer's
theorem~\cite{MR96a:55019} that their $B\Z/p$-nullification is not
contractible. We prove in Proposition~\ref{prop covers} that,
under this finiteness condition, there are basically no other
$H$-spaces with infinitely many non-trivial homotopy groups than
the highly connected covers of mod $p$ finite $H$-spaces.

\vskip 3mm

\noindent {\bf Acknowledgements.} We would like to thank Bill
Dwyer and Clarence Wilkerson for attracting our attention to this
problem.

\section{The homotopy groups of $H$-spaces}
\label{sec H}
The original theorem \cite[Theorem~1.3]{MR92b:55004} by Dwyer and
Wilkerson about the homotopy groups of certain $2$-connected
spaces relies on the equivalence between a cohomological condition
and a topological one. Namely, the loop space of a $p$-complete
space is $B\Z/p$-null if and only if the module of indecomposable
elements in mod $p$ cohomology is locally finite,
\cite[Proposition~3.2]{MR92b:55004}. In fact, this result can be
understood as a reduction step to the theorem of Lannes and
Schwartz: If $X$ is $2$-connected and $QH^*(X; \F_p)$ is locally
finite, then $\Omega X$ is $B\Z/p$-null as we just recalled; thus
the cohomology $H^*(\Omega X; \F_p)$ is locally finite, which
implies by ~\cite{MR827370} that $\Omega X$ has infinitely many
non-trivial homotopy groups (unless it is contractible).

When $X$ is an $H$-space, we were able to obtain an extension of
\cite[Proposition~3.2]{MR92b:55004} using the Krull filtration of
the category of unstable modules over $\A_p$.

\begin{Thm} \cite[Theorem~5.3]{CCS}
\label{thm CCS}
Let $X$ be a connected $H$-space such that $T_V H^* (X; \F_p)$ is
of finite type for any elementary abelian $p$-group $V$. Then
$QH^*(X; \F_p) \in \U_n$ if and only if $\Omega^{n+1} X$ is a
$B\Z/p$-null space. \hfill{\qed}
\end{Thm}

We obtain then, as in~\cite{MR92b:55004}, a result on the homotopy
groups of sufficiently connected spaces satisfying the conditions
of our theorem.

\begin{Thm}
\label{prop Hdichotomy}
Let $X$ be an $(n+2)$-connected $H$-space for some integer $n \geq
0$ such that $T_V H^* (X; \F_p)$ is of finite type for any
elementary abelian $p$-group $V$. Assume that $QH^*(X; \F_p)$ lies
in $\U_n$. Then either $X$ is contractible, or it has infinitely
many non-trivial homotopy groups with $p$-torsion. In the second
case the iterated loop space $\Omega^{n+1} X$ has infinitely many
non-trivial $k$-invariants.
\end{Thm}

\begin{proof}
Assume that $X$ is not contractible. By Theorem~\ref{thm CCS},
$\Omega^{n+1} X$ is a $B\Z/p$-null space. We know thus from
\cite[Theorem~5.5]{CCS} (compare with \cite[Theorem~7.2]{B2}) that
the homotopy fiber $F$ of the nullification map $X \rightarrow
P_{B\Z/p} X$ is a $p$-torsion Postnikov piece with its homotopy
groups concentrated in degrees from $1$ to $n+1$.

We infer from the homotopy long exact sequence that $P_{B\Z/p} X$
is simply connected and not contractible, because $X$ is
$(n+2)$-connected and not contractible. Therefore, the
Lannes-Schwartz theorem~\cite{MR827370} applies. The space
$P_{B\Z/p} X$ must have an infinite number of non-trivial homotopy
groups with $p$-torsion, and so does $X$.

The assertion about the $k$-invariants follows from the fact that
an Eilenberg-Mac Lane space $K(A, m)$ is not $B\Z/p$-local if $A$
contains $p$-torsion.
\end{proof}

\begin{Cor}\label{coversoflocals}
Let $n\geq 0$ and $X$ be a $p$-complete $H$-space such that
$T_VH^*(X; \F_p)$ is of finite type, $H^*(X; \F_p)$ is
$(n+2)$-connected and $QH^*(X; \F_p)\in \U_{n}$. Then $X$ is the
$(n+2)$-connected cover of a $B\Z/p$-null $H$-space. \hfill{\qed}
\end{Cor}

\begin{proof}
As $X$ is $p$-complete, the connectivity condition on the mod $p$
cohomology implies that $X$ itself is $(n+2)$-connected. The
fibration $F \rightarrow X \rightarrow P_{B\Z/p} X$ used in the
proof of Theorem~\ref{prop Hdichotomy} exhibits now $X$ as a
highly connected cover of a $B\Z/p$-null space.
\end{proof}

\section{On $H$-spaces that are Postnikov pieces}
\label{sec fgoverAp}

Whereas algebraic conditions that ensure that a space is not a
Postnikov piece are frequently encountered in the literature, a
characterization of Postnikov pieces in terms of their cohomology
seems out of reach. Our aim in this section is to provide a
satisfactory answer for $H$-spaces. Let us first look at the
cohomology of an $H$-space with finitely many $p$-torsion homotopy
groups.

\begin{Prop}
\label{prop Postnikovpieces}
Let $X$ be an $H$-space which is a $p$-torsion Postnikov piece of
finite type. Then $H^* (X; \F_p)$ is a finitely generated algebra
over~$\A_p$.
\end{Prop}

\begin{proof}
The homotopy group of a $p$-torsion Eilenberg-Mac Lane space of
finite type is a finite direct sum of cyclic groups $\Z/p^n$ and
Pr\"ufer groups $\Z_{p^\infty}$. The cohomology of such spaces has
been computed by Cartan and Serre. It is finitely generated as an
algebra over~$\A_p$ (see for example
\cite[Section~8.4]{MR95d:55017}).

In~\cite[Proposition~6.2]{CCS} we proved that the cohomology of
the total space of an $H$-fibration is a finitely generated
algebra over $\A_p$, if so are the cohomology of both the fiber
and the base. Therefore the result follows by induction on the
number of homotopy groups of the $H$-space~$X$.
\end{proof}

We offer now our characterization by combining
Proposition~\ref{prop Postnikovpieces} with a result analogous to
\cite[Lemma~2.1]{MR96a:55019} on the $B\Z/p$-nullifica\-tion of
$p$-torsion Postnikov pieces.

\begin{Thm}
\label{thm Postnikovpieces}
Let $X$ be an $H$-space. Then $X$ is a $p$-torsion Postnikov piece
of finite type if and only if $P_{B\Z/p} X$ is contractible and
$H^* (X; \F_p)$ is a finitely generated algebra over~$\A_p$.
\end{Thm}

\begin{proof}
If $H^*(X; \F_p)$ is a finitely generated algebra over $\A_p$,
then by~\cite[Lemma~6.1]{CCS} the module $QH^*(X; \F_p)$ belongs
to $\U_n$ for some $n$ and $T_V H^*(X; \F_p)$ is of finite type
for any~$V$. Therefore Theorem~\ref{thm CCS} applies, so
$\Omega^{n+1} X$ is $B\Z/p$-null. Now using Bousfield's
description~\cite[Theorem~7.2]{B2} of the homotopy fiber of the
nullification map $X \rightarrow P_{B\Z/p} X$ (see
~\cite[Theorem~5.5]{CCS} in this concrete setting), this fiber is
a $p$-torsion Postnikov piece. As $P_{B\Z/p} X$ is contractible,
$X$ itself is a Postnikov piece, and since $H^*(X; \F_p)$ is
finitely generated, an elementary Serre spectral sequence argument
shows that $X$ is of finite type.

Conversely, if $X$ is a $p$-torsion Postnikov piece which is an
$H$-space, then its $B\Z/p$-nullification is contractible. This
statement follows from the fact that $p$-torsion Eilenberg-Mac
Lane spaces are $B\Z/p$-acyclic and $P_{B\Z/p}$ preserves
fibrations in which the fiber is $B\Z/p$-acyclic
(\cite[Theorem~1.H.1]{Dror}). We conclude by Proposition~\ref{prop
Postnikovpieces}.
\end{proof}

\begin{Rmk}
\label{rmk BS3}
When $X$ is not an $H$-space, this characterization fails.
Consider for example $X$, the homotopy fiber of the nullification
map $BS^3\rightarrow P_{B\Z/p} BS^3 \simeq \Z[1/p]_\infty BS^3$
(see \cite[Theorem~1.7, Lemma~6.2]{MR97i:55028}). Then $P_{B\Z/p}
X$ is contractible by \cite[Theorem~1.H.2]{Dror}, and $H^*(X;
\F_p)$ is isomorphic to $H^*(BS^3; \F_p)$, hence finitely
generated as an algebra. Notwithstanding $X$ is not a $p$-torsion
Postnikov piece (see also \cite[Example 3.7]{CCS}).
\end{Rmk}

We wish to mention that there is an obvious way to apply our
results to spaces that are not $H$-spaces, namely by considering
their loop space.

\begin{Cor}
\label{cor covers}
A $1$-connected space $X$ is a $p$-torsion Postnikov piece of
finite type if and only if $P_{\Sigma B\Z/p} X$ is contractible
and $H^*(\Omega X; \F_p)$ is finitely generated as an algebra
over~$\A_p$.
\end{Cor}

\begin{proof}
A space is a Postnikov piece if and only if its loop space is so.
Thus Theorem~\ref{thm Postnikovpieces} applies to the connected
$H$-space $\Omega X$ and we conclude since $P_{B\Z/p} \Omega X
\simeq \Omega P_{\Sigma B\Z/p} X$, \cite[Theorem~3.A.1]{Dror}.
\end{proof}

It would be nice to find a characterization in terms of the
cohomology of $X$ itself rather than the mod $p$ loop space
cohomology.

\section{Connected covers of finite $H$-spaces}

This section is devoted to analyze the nature of $H$-spaces whose
mod $p$ cohomology is finitely generated over $\A_p$ but that are
not Postnikov pieces. Examples of $H$-spaces having finitely
generated cohomology as an algebra over the Steenrod algebra are
the highly connected covers of simply connected mod~$p$ finite
$H$-spaces (such as odd dimensional spheres completed at odd
primes). Such spaces have obviously infinitely many non-trivial
homotopy groups, as a direct consequence of Serre's original
theorem \cite[Th\'eor\`eme~10]{MR0060234} and its generalization
given by McGibbon and Neisendorfer \cite[Theorem~1]{MR749108}. We
prove that there are basically no other $H$-spaces which do have
infinitely many non-trivial homotopy groups: Any $H$-space with
finitely generated mod~$p$ cohomology as an algebra over the
Steenrod algebra differs from a mod~$p$ finite one by only a
finite number of homotopy groups. In other words, some iterated
loop space of such an $H$-space coincides with the iterated loop
space of a mod~$p$ finite $H$-space.

\begin{Prop}
\label{prop covers}
Let $X$ be an $H$-space such that $H^* (X; \F_p)$ is a finitely
generated algebra over $\A_p$. Then there exist an integer $n$ and
an $H$-space $Y$ with finite mod $p$ cohomology such that the
$(n+2)$-connected cover of $Y$ and $X$ are equivalent. Moreover,
when $P_{B\Z/p}X$ is not contractible, $X$ has infinitely many
non-trivial homotopy groups.
\end{Prop}

\begin{proof}
The space $Y$ is obtained as the $B\Z/p$-nullification of $X$. As
we explain in \cite[Proposition~6.8]{CCS} its mod $p$ cohomology
is finite because it is both finitely generated as an algebra over
$\A_p$ and locally finite as an unstable module. The integer $n$
is the smallest one such that $QH^* (X; \F_p)$ belongs to $\U_n$,
which exists since $H^* (X; \F_p)$ is a finitely generated algebra
over $\A_p$. Moreover $X$ and $Y$ differ by a finite number of
homotopy groups (concentrated in degrees from $1$ to $n+1$)
because the homotopy fiber of $X \rightarrow Y$ is a $p$-torsion
Postnikov piece, see \cite[Theorem~5.5]{CCS}.
\end{proof}

\section{A variation with Neisendorfer's functor}
\label{neisendorfer}

In Section~\ref{sec fgoverAp} we considered the nullification
functor $P_{B\Z/p}$. Next we explain how to obtain analogous
results for the functor $(P_{B\Z/p} )^\wedge_p$ introduced by
Neisendorfer in~\cite{MR96a:55019}. However we need to add an
extra condition on the fundamental group because $S^1$ is a
$B\Z/p$-null space.

\begin{Prop}
\label{rmk Neisendorfer}
Let $X$ be a $p$-complete $H$-space with finite fundamental group.
Then $X$ is a Postnikov piece of finite type if and only if
$(P_{B\Z/p} X)^\wedge_p$ is contractible and $H^* (X; \F_p)$ is a
finitely generated algebra over~$\A_p$.
\end{Prop}

\begin{proof}
If $H^* (X; \F_p)$ is a finitely generated algebra over~$\A_p$,
consider the fibration $F \rightarrow X \rightarrow P_{B\Z/p} X$.
Since $H$-fibrations are preserved by $p$-completion and
$(P_{B\Z/p} X)^\wedge_p$ is contractible, we see that $X\pcom$ is
a $p$-complete Postnikov piece (more precisely, the $p$-completion
of a $p$-torsion Postnikov piece).

Conversely, if $X$ is a connected $p$-complete Postnikov piece of
finite type with finite fundamental group, then
$(P_{B\Z/p}X)^\wedge_p$ is contractible by Neisendorfer's
result~\cite[Lemma~2.1]{MR96a:55019} and the cohomology is
finitely generated as an algebra over~$\A_p$ by Proposition
\ref{prop Postnikovpieces}.
\end{proof}

The connectivity assumption in Theorem~\ref{prop Hdichotomy}
cannot be relaxed because of the obvious example of $K(\Z, n+2)$.
In fact this is essentially the unique $(n+1)$-connected $H$-space
which is a Postnikov piece such that $QH^*(X; \F_p)$ lies in
$\U_n$.

\begin{Prop}
\label{prop KZn}
Let $X$ be an $(n+1)$-connected $H$-space for some integer $n \geq
0$ such that $T_V H^* (X; \F_p)$ is of finite type for any
elementary abelian $p$-group~$V$. Assume that the module of
indecomposable elements $QH^*(X; \F_p)$ lies in $\U_n$ and that
$X$ is a Postnikov piece. Then $X$ is, up to $p$-completion,
homotopy equivalent to the product of finitely many copies of
$K(\Z^\wedge_p, n+2)$.
\end{Prop}

\begin{proof}
Consider the fibration $F \rightarrow X \rightarrow P_{B\Z/p} X$
as in the proof of Theorem~\ref{prop Hdichotomy}. Since
$\Omega^{n+1} X$ is a $B\Z/p$-null space, we know that the fiber
$F$ is a $p$-torsion Postnikov piece and its homotopy groups are
concentrated in degrees from $1$ to $n+1$. By Proposition~\ref{rmk
Neisendorfer} $(P_{B\Z/p} X)^\wedge_p$ is contractible, so $X$
itself is, up to $p$-completion, homotopy equivalent to~$F$. The
connectivity assumption implies that $F^\wedge_p$ must be
$(n+1)$-connected. Thus the only non-trivial homotopy group of $F$
is $\pi_{n+1} F$ and it must be a finite product of copies of
$\Z_{p^\infty}$, since $K(\Z_{p^\infty}, n+1)^\wedge_p \simeq
K(\Z^\wedge_p, n+2)$.
\end{proof}

Finally we propose a characterization of Postnikov pieces which
are infinite loop spaces.

\begin{Prop}
\label{prop McGibbonNeisendorfer}
Let $X$ be a $p$-complete infinite loop space with finite
fundamental group. Then $X$ is a Postnikov piece of finite type if
and only if $H^* (X; \F_p)$ is a finitely generated algebra
over~$\A_p$.
\end{Prop}

\begin{proof}
The $B\Z/p$-nullification of a connected infinite loop space with
$p$-torsion fundamental group is trivial up to $p$-completion by
McGibbon's main theorem in~\cite{MR97g:55013}. We conclude by
Proposition~\ref{rmk Neisendorfer}.
\end{proof}


\providecommand{\bysame}{\leavevmode\hbox
to3em{\hrulefill}\thinspace}
\providecommand{\MR}{\relax\ifhmode\unskip\space\fi MR }
\providecommand{\MRhref}[2]{%
  \href{http://www.ams.org/mathscinet-getitem?mr=#1}{#2}
} \providecommand{\href}[2]{#2}



\bigskip
{\small
\begin{center}
\begin{minipage}[t]{8 cm}
Nat\`{a}lia Castellana and Jer\^{o}me Scherer\\ Departament de Matem\`atiques,\\
Universitat Aut\`onoma de Barcelona,\\ E-08193 Bellaterra, Spain
\\\textit{E-mail:}\texttt{\,Natalia@mat.uab.es}, \\
\phantom{\textit{E-mail:}}\texttt{\,\,jscherer@mat.uab.es}
\end{minipage}
\begin{minipage}[t]{7 cm}
Juan A. Crespo \\ Departament de Economia i de Hist\`{o}ria
Econ\`{o}mica,
\\ Universitat Aut\`onoma de Barcelona,\\ E-08193 Bellaterra,
Spain
\\\textit{E-mail:}\texttt{\,JuanAlfonso.Crespo@uab.es},
\end{minipage}
\end{center}}

\end{document}